\documentclass{article}
\usepackage{amsthm,amsfonts,amssymb}
\usepackage[utf8]{inputenc}
\usepackage[T1]{fontenc}
\input amssym.def
\input amssym.tex
\def\newline{\hfill\break}
\def\scong{{\scriptstyle\|}\lower.2ex\hbox{$\wr$}}
\def\Z{{\Bbb Z}}
\def\R{{\Bbb R}}
\def\Q{{\Bbb Q}}

\def\C{{\Bbb C}}

\def\End{\mathop{\rm End}\nolimits}

\def\Br{\mathop{\rm Br}\nolimits}

\def\Cor{\mathop{\rm Cor}\nolimits}

\title{Mixed Characteristic Artin Schreier Polynomials}
\author{David J. Saltman\\  
Center for Communications Research\\
805 Bunn Drive\\
Princeton, NJ 08540}

\newtheorem{theorem}{Theorem}[section]
\newtheorem{corollary}[theorem]{Corollary}
\newtheorem{lemma}[theorem]{Lemma}
\newtheorem{proposition}[theorem]{Proposition}

\begin{document}

\maketitle

\begin{abstract}
We present here a version of the Artin-Schreier polynomial that works in 
any characteristic. Let $C_p$ be the cyclic group of prime order $p$. 
Equivalently, we prove one can lift degree $p$ 
cyclic $C_p$ extensions over local rings $R,M$ where $R/M$ has characteristic 
$p$ and $R$ has arbitrary characterstic. 

Let $\rho \in \C$ be a primitive $p$ root of one. We first consider the case $R$ has a 
primitive $p$ root of one, by which we mean that there is a given homomorphism 
$f: \Z[\rho] \to R$. In this context we can write a specific polynomial 
which generalizes the Artin-Schreier polynomial. We next consider arbtitrary 
$R$ and construct a mixed characteristic ``generic'' Galois extension that proves the 
lifting result, but here we do not supply a polynomial. 

It is useful to view these results in terms of Galois actions. 
If $C_p$ is generated by $\sigma$, then Artin-Schreier polynomials 
describe $C_p$ Galois extensions generated by elements $\theta$ where 
$\sigma(\theta) = \theta + 1$. If $\rho' = f(\rho)$, then our lifted Galois 
action is $\sigma(\theta) = \rho'\theta + 1$. In section two we descend 
this extension to rings without root of one. In section three we 
use this new description of cyclic extensions, give the obvious 
description of the corresponding cyclic algebras, and then define 
a new algebra which specializes to general differential crossed products 
in characteristic $p$ but which are cyclic in characteristic not $p$. 
The algebra in question is generated by $x,y$ subject to the relation 
$xy - \rho{yx} = 1$. 
\end{abstract}

\section*{Introduction} 

The Artin-Schreier polynomial $x^p - x - a$ 
is deservedly famous. If $F$ is a field of 
characteristic $p$, then any Galois extension 
$L/F$ of degree $p$ is described by such a 
polynomial for an appropriate choice of $a \in F$. 
Moreover, the action of the Galois group $G$ of this 
extension is particularly nice. Namely,  
$G$ has a generator $\sigma$ such that if $\theta$ 
is a root of this polynomial, then $\sigma(\theta) = \theta + 1$. 

It is less well known that this class of polynomials is even stronger. Suppose $R$ 
is a commutative ring of characteristic $p > 0$, 
and $S/R$ is cyclic Galois of degree $p$, 
then $S \cong R[x]/(x^p - x - a)$. 
That is, the Artin Schreier polynomials more 
generally
describes Galois extensions over commutative 
rings. Note that the key fact here is 
that $x^p - x - a$ is separable for all $a$. 
Compare this with $x^n - a$, in characteristic 
prime to $n$, which is only separable when 
$a \not= 0$. This is not an accident 
of the choice of polynomial. The polynomial 
ring $F[a]$, for $F$ of characteristic prime to 
$n$, has no Galois extensions of degree $n$ 
that do not come from $F$. 

Of course, these polynomials are one example 
of so called generic polynomials and define so called generic Galois extensions. Recall 
that if there is a generic Galois extension 
for a group $G$ over $F$, $R$ is a local $F$ 
algebra with residue field $R/M = K$, 
and $L/K$ is a $G$ Galois extension, then there is 
a $G$ Galois $S/R$ such that $S/MS \cong L$. 
(In fact all we need know about $R \to R/M$ 
is that it is surjective on units.)
This is called the lifting property. 
The stronger properties of the Artin Schreier 
polynomial imply a stronger lifting property. 
Namely, if $F$ has characteristic $p$ 
and $\phi: R \to S$ is a surjection of commutative 
$F$ algebras, $G = C_p$ is the cyclic group 
of order $p$, and $S'/S$ is $C_p$ Galois, 
then there is a $C_p$ Galois $R'/R$ with 
$S' \cong R' \otimes_{\phi} S$. 

With the above in mind, it is natural to ask 
about lifting $C_p$ Galois extensions in 
mixed characteristics or where $R$ does 
not contain a field. This question came up quite 
concretely in the paper \cite{Su} by Suresh, and 
this paper starts from some of his observations. 

The change of variables that is the key to the first 
section has been used before. It was central to 
the arguments in \cite{Su}, and also appeared in the 
work on the Oort conjecture (\cite {OSS}) 
reviewed in \cite {O}. The explicit Galois action, 
the generality of the description, 
and the removal of dependence on a root of one 
in section two seems new. Since we have an explicit Galois 
action, namely $\theta \to \rho\theta + 1$, that gives an alternate description of cyclic algebras 
of degree $p$. However, we get the best result by following the 
lead from differential crossed product algebras in characteristic 
$p$, and defining a "cyclic" maximal order, characteristic free, 
with very good ramification properties. The key relationship 
is $xy - \rho{yx} = 1$. This defines the "curious" algebra 
discussed in section three.  

To put our results in context, if $G$ is any $p$ group then over characteristic 
$p$ fields there is a generic $G$ Galois extension \cite{S}. However for arbitrary $G$ 
these are known to not exist in characteristic 0 and so our mixed characterstic 
generic extensions will certainly not exist for general $p$ groups $G$. 
However, there is a large range of groups where they might exist, the obvious first 
case being cyclic groups of order $p^n$. This is the subject of future investigation. 

\section{The polynomial} 

Let $\rho$ be the primitive $p$ root of one 
over $\Z$. That is, $\Z[\rho] = 
\Z[x]/(x^{p-1} + x^{p-2} + \cdots + 1)$. 
Of course $\Z[\rho]$ is a Dedekind domain 
containing $\eta = \rho - 1$ and $\eta$ 
is a prime totally ramified over $p$. 

We need to make this explicit. 
$$1 = (\eta + 1)^p  = \eta^p + p\eta\left(\sum_{i=1}^{p-1} b_i\eta^{i-1}\right) + 1 = \eta^p + p\eta{y} + 1
$$ 
where $b_i$ is an integer and the appropriate binomial coefficient 
divided by $p$ and $y = \sum_i b_i\eta^{i-1} 
\in \Z[\rho]$. Since $b_1 = 1$, 
$y \in 1 + \eta\Z[\rho]$. Thus 
$\eta^p + p\eta{y} = 0$ or $\eta^{p-1} = -py$  
or $p = x\eta^{p-1}$ where $x = -1/y$. 
Modulo $p$, $y$ and hence $x$ are units. 
If $q \not= p$, then modulo any prime containing 
$q$, $p$ and $\eta$ are units and so $y$ is a unit. 
It follows that $x$ is a unit in $\Z[\rho]$. 
Since $y$ is congruent to 1 modulo $\eta$, 
$x$ is congruent to $-1$ modulo $\eta$. 

Now we consider the polynomial ring $R_0 = 
\Z[\rho][u]$ and the polynomial 
$X^p - (1 + u\eta^p)$. For reasons that will soon be clear, 
we set $R = R_0(1/(1 + u\eta^p))$. If we substitute 
$X = 1 + \eta{Z}$ we get 

$$(1 + \eta{Z})^p - (1 + u\eta^p) = 
\eta^p{Z^p} - u\eta^p + 
p\eta\left(\sum_1^{p-1} b_iZ^i\eta^{i-1}\right)$$ 
$$= \eta^p{Z^p} - u\eta^p + 
x\eta^p\left(\sum_1^{p-1} b_iZ^i\eta^{i-1}\right) = 
\eta^p\left(Z^p + x\sum_1^{p-1} b_iZ^i\eta^{i-1} - u\right)$$
and we define $g(Z)$ so that this last expression 
is $$\eta^p(Z^p + g(Z) - u)$$ 
and we have also shown that 
$$(1 + z\eta)^p = 1 + (g(z) + z^p)\eta^{p}$$
Since $b_1 = 1$ we note that $g(Z)$ is congruent
to $-Z$ modulo $\eta$. That is, our new 
polynomial $Z^p + g(Z) - u$ is the Artin-Schreier 
polynomial $Z^p - Z - u$ modulo $\eta$. 

Suppose $\alpha$ is the image of $X$ 
in $S = R_0[X]/(X^p - (1 + u\eta^p))$. 
Then $\theta = (\alpha - 1)/\eta \in S(1/\eta)$ is a root of $Z^p + g(Z) - u$. Of course, 
$Z^p + g(Z) - u$ has coefficients in $R_0$ and 
so $\theta$ is integral over $R_0$. 
Let $T_0 = R_0[Z]/(Z^p + g(Z) - u)$ and view 
$S_0 \subset T_0$ by mapping $X$ to $1 + Z\eta$. 
Since $\Z[\rho]/(\eta) = F_p$, the field of 
$p$ elements, $T_0/\eta{T_0} = F_p[u][Z]/(Z^p - Z - u)$ 
and so $\eta$ is a prime element of $T_0$. 

Let $S = S_0(1/(1 + u\eta^p)) = R[X]/(X^p - (1 + u\eta^p))$ and 
$T = T_0(1/(1 + u\eta^p)) = R[Z]/(Z^p + g(Z) - u)$. 
Define $\sigma: S \cong S$ by letting 
$\sigma$ be the identity on $R$ and setting 
$\sigma(\alpha) = \rho\alpha$. Note that $\sigma(S_0) = S_0$. 
Then 
over $S$, $X^p - (1 + u\eta^p) = 
\prod_i (X - \sigma^i(\alpha)) = 
\prod_i (X - \rho^i\alpha)$. 
Under our identification, $T \subset S(1/\eta)$.  
The map 
$\sigma$ extends to $S(1/\eta)$ and 
$\sigma(\theta) = (\rho(\alpha) - 1)/\eta = 
((\eta + 1)\alpha - 1)/\eta = \alpha + \theta = 
\eta\theta + 1 + \theta = \rho\theta + 1$. 
That is, 

$$\sigma(\theta) = \rho\theta + 1.$$ 

Thus $\sigma(T_0) = T_0$ and $\sigma(T) = T$. 
Of course, modulo $\eta$, this is just the 
Artin-Schreier action.  
By induction $\sigma^i(\theta) = \rho^i\theta + 
\rho^{i-1} + \cdots + 1$. Also, 
$Z^p + g(Z) - u = \prod_i (Z - \sigma^i(\theta))$.

\begin{proposition}\label{prop:1.1} Suppose $f: \Z[\rho] \to R'$ defines 
$R'$ as a $\Z[\rho]$ algebra. Let $g: R \to R'$ be any 
$\Z[\rho]$ algebra homomorphism and set $a = f(u)$. 
Then $g(Z^p + g(Z) - u) = Z^p + g(Z) - a$ is separable if and only 
if $1 + a\eta^p$ is invertible. In particular, the extension $T/R$ above 
is Galois with group generated by $\sigma$. 
\end{proposition}

\begin{proof} Let $\theta$ be a root of $Z^p + g(Z) - u$ 
and $\Delta$ the discriminant. Then $\Delta$ 
the product of expressions of the form 
$\sigma^j(\sigma^i\theta - \theta) = 
\sigma^j((\rho^i - 1)\theta + \rho^{i-1} + \cdots + 1) = (\rho^{i-1} + \cdots + 1)(\eta\theta + 1)$ 
where $i \not= 0$ and $j = 0, \ldots, p-1$. 
Fixing $i$ and taking the product for all $j$ 
yields $(\rho^{i-1} + \cdots + 1)(1 + u\eta^p)$ 
since $\eta\theta + 1$ is a root of 
$X^p - (1 + u\eta^p)$. Also, 
$(\rho^{i-1} + \cdots + 1)$ is $i$ modulo $\eta$ 
and in any finite field of characteristic not $p$ 
this equals $(\rho^i - 1)/(\rho - 1) \not= 0$. 
Thus $(\rho^{i-1} + \cdots + 1)$ is a unit in 
$\Z[\rho]$ proving the first statement. 

As for the second statement, we use the 
criterion of \cite[page~81]{DI}. Since 
$\sigma^i(\theta) - \theta$ is a unit, 
we only need to verify that the fixed ring $T^{\sigma} = 
R$. If $K$ is the field of fractions 
of $R$ and $L = T \otimes_{R} K$ 
then $L^{\sigma} = K$ because $L/K$ is the 
Kummer field extension $K((1 + u\eta^p)^{1/p})/K$. 
Since $T/R$ is integral, and 
$R$ is integrally closed, $T^{\sigma} = 
R$ is clear.
\end{proof}

In the generic case $\eta$ is a non--zero divisor and one can pass easily from 
$\theta$ to $\alpha = \eta\theta + 1$. 
In general we cannot, and so we must develop 
relations in the generic case and then 
specialize them. Furthermore, the notation 
we develop will be quite suggestive. In the following discussion 
we will operate in two domains. If we work ``generically'' we will 
be making our definitions and calculations in $R_0$ or $R$ or an extension 
ring $R[x,y]$ etc. The key point will be that $\eta$ is a 
non--zero divisor. When we work generally, we mean we have an arbitrary 
$\Z[\rho]$ algebra $R'$ and a homomorphism $g:R_0 \to R'$ or $R[x,y]\to R'$ etc. 
 
Let $A$ be the monoid $1 +\eta{R}$ under 
multipication, and $A^*$ the subgroup of invertible elements.  
Define $\phi: R \to A$ via $\phi(x) = 
1 + x\eta$. Since $\eta$ is a non--zero divisor, 
this is a bijection.  
To make this some sort of homomorphism we define 
$x \oplus y$ by writing $(1 + x\eta)(1 + y\eta) = 
 1 + (x \oplus y)\eta$ or $\phi(x \oplus y) = 
\phi(x)\phi(y)$. 
Note that we are just using helpful notation. 
Clearly, $x \oplus y = x + y + xy\eta$. 
If we look at general $g: R \to R'$ (which includes $R \to R[x,y]$) 
we define 
$\eta' = f(\eta)$, $A' = 1 + \eta'R'$, $\phi: R' \to A'$, 
$\phi(x \oplus y) = \phi(x)\phi(y)$ etc. In the following 
if we make an observation about $R$ or $R[x,y]$, we will not restate the obvious 
inference about $R'$. 

We can define $x \ominus y = z$ to mean $x = y \oplus z$. 
Of course $z$ 
exists if $y \in A'^*$. Note that since 
$\rho' = 1 + \eta'$, $\phi(1) = \rho'$. 
Also, the basic relation between the Kummer 
element $\alpha$ and $\theta$ described above 
is just that $\phi(\theta) = \alpha$. 
Thus $\sigma(\alpha) = \rho\alpha$ translates 
into $\sigma(\phi(\theta)) = \phi(1)\phi(\alpha)$ 
or, suggestively, 
$$\sigma(\theta) = 1 \oplus \theta$$ 
which makes sense because $1 \oplus \theta = 
1 + \theta + \theta\eta = 1 + \rho\theta$. 

We are interested in the case $Z^p + g(Z) - a$ defines a 
Galois extension so let $g: R \to R'$ be as above. 
Let $S' = R'[Z]/(Z^p + g(Z) - a)$ and let $\theta$ be the image of 
$Z$, so $S'/R'$ is $C_p$ Galois and $\sigma(\theta) = \theta \oplus 1$. 
The next result describes how we can change 
the choice of $\theta$ and then track the corresponding 
change in $a$. To describe the result, 
let $A_p$ be the monoid $1 + R'\eta'^p$ 
and $\phi_p: R' \to A_p$ the map 
$r \to 1 + r\eta^p$. In a similar way as before define $x \oplus_p y$ by 
$\phi_p(x)\phi_p(y) = \phi_p(x \oplus_p y)$. 
The point of this and subsequent definitions 
is a to distinguish $z$ from it's image 
$1 + z\eta$ or $1 + z\eta^p$. 
This is crucial when $\eta$ is a zero 
divisor. 
We gathered some easy manipulation rules for 
these operations in the appendix. 

\begin{proposition}\label{prop:1.2}  If $\theta$ is as 
above and $\theta'$ also satisfies 
$\sigma(\theta') = 1 \oplus \theta'$ then 
$\theta' = \theta \oplus z$ for some $\sigma$ 
fixed $z$. If $\theta$ satisfies $Z^p + g(Z) - a = 0$ then $\theta'$ satisfies 
$Z^p + g(Z) - a' = 0$ 
where $a' = a \oplus_p (g(z) + z^p)$. 
Conversely, if $a' = a \oplus_p (g(z) + z^p)$ for some $z$ with $\phi(z) \in A^*$ then 
$Z^p + g(Z) - a$ and $Z^p + g(Z) - a'$ define
isomorphic $C_p$ Galois extensions. 
\end{proposition}

\begin{proof} $(1 + \eta\theta)^p = 1 + a\eta^p$ and so 
$1 + \eta\theta$ is invertible.  We can 
set $z = \theta' \ominus \theta$ and the first 
statement is clear. To show the second statement, 
we can work in the generic case and hence 
assume $\eta$ is a non--zero divisor. We compute that 
$(1 + \theta'\eta)^p = 
(1 + \theta\eta)^p(1 + z\eta)^p = 
(1 + \theta\eta)^p(1 + \eta^p(g(z) + z^p))^p = 
(1 + u\eta^p)(1 + (g(z) + z^p)\eta^p) = 1 + (u \oplus_p (g(z) + z^p))\eta^p)$. 
Thus $\theta'$ behaves as claimed. 
For the last statement we observe that 
$\tau(\theta') = \theta \oplus z$ defines the 
isomorphism.
\end{proof} 

For $s$ a positive integer define $s * a$ via $\phi(a)^s = \phi(s * a)$. That is, $s * a = 
a \oplus \cdots \oplus a$ ($s$ times). 
It follows that $s * (a \oplus b) = (s * a) \oplus (s * b)$. 
Note that $\phi(s*1) = \rho^s = 1 + \delta_s\eta$ where $\delta_s = 1 + \cdots + \rho^{s-1}$. 
That is, $$s*1 = \delta_s.$$ 
Now 
$\sigma(s*\theta) = s*\sigma(\theta) = 
s*(1 \oplus \theta) = (s*1) \oplus (s*\theta) = 
\delta_s \oplus (s*\theta)$ so 
$$\sigma(s*\theta) = s*\theta \oplus \delta_s.$$

The operation $p*z$ is not so useful because, as we see next, 
$\eta^{p-1}$ divides $p*z$. 
As mentioned above, this is problematic when 
$\eta$ is a zero divisor. 
In fact, 
$(1 + z\eta)^p = 1 + p\eta{y} + z^p\eta^p$ 
and since $p = x\eta^{p-1}$ we have that 
$(1 + z\eta)^p \in A_p$. Working in the generic case, 
we can define 
$p \#_p z$ by setting $(1 + z\eta)^p = 
1 + (p \#_p z)\eta^p$. Thus $p*z = (p\#_p z)\eta^{p-1}$. More generally if 
$pr$ is a multiple of $p$ 
we define $pr \#_p z$ via $(1 + (pr \#_p z)\eta^p) = (1 + z\eta)^{pr}$. We can restate 
our basic equation: 

\begin{lemma}\label{lem:1.3}  If $\sigma(\theta) = 
\theta \oplus 1$ then $p \#_p \theta = g(\theta) + 
\theta^p$ and is $\sigma$ fixed. 
\end{lemma}

\begin{proof} We can prove this in the generic case 
and hence where $\eta$ is a non zero divisor. 
Now all this says (as we saw above) that 
$(1 + \eta\theta)^p$ is $\sigma$ fixed and 
equals $1 + u\eta^p$ where $u = \theta^p + 
g(\theta)$.
\end{proof}

We can easily derive more rules for dealing 
with these new operations and we have also gathered them in the 
appendix. 

Note that $\delta_s = \delta_{s + rp}$ 
so we can view the subscript $s$ of $\delta_s$ 
as an element of $\Z/p\Z$. 

The next results shows the extent to which 
$T/R$ above is a generic object. 
Since our extension subsumes the prime to $p$ 
characteristic and hence Kummer case, we 
must restrict the base ring $R'$ and also 
need $R'$ to ``contain'' $\rho$, by which we 
mean that $R'$ is a $\Z[\rho]$ algebra. 
That is, there is a fixed ring homomorphism 
$f: \Z[\rho] \to R'$.  

\begin{theorem}\label{thm:1.4} Let $R'$ be a semilocal 
commutative $\Z[\rho]$ algebra defined 
by $f:\Z[\rho] \to R'$. Let $S'/R'$ 
be a $C_p$ Galois extension. Then there is 
an extension $g: R \to R'$ such that 
$S' \cong T \otimes_{\phi} R'$ as 
$C_p$ Galois extensions. 
\end{theorem}

\begin{proof} Recall $\rho' = f(\rho)$. Our assumptions on $R'$ 
imply that $S'/R'$ has a normal basis. 
That is, $S' \cong R'[C_p]$ as a $R'[C_p]$ 
module. If $z,\sigma(z),\ldots, 
\sigma^{p-1}(z)$ is this normal basis, 
note that $a = \sum_i \sigma^i(z)$ must be invertible and so multiplying by $1/a$ 
we may assume $\sum_i \sigma^i(z) = 1$. 

Set $r_{p-i} = f(\delta_i)$, 
so $r_0 = 1$, $r_1 = 0$ and $r_i = 
\rho'{r_{i+1}} + 1$. If we set $\theta' = 
\sum_{i=0}^{p-1} r_i\sigma^i(z)$, 
then $\sigma(\theta') = \sum_i r_i\sigma^{i+1}(z) = 
\rho'\theta' + 1$. It is clear that $1,\theta',
\sigma^2(z),\ldots,\sigma^{p-1}(z)$ form a basis 
of $S'/R'$ but we need that $1,\theta',\theta'^2,\ldots, \theta'^{p-1}$ is a basis. 
\end{proof}

\begin{lemma}\label{lem:1.5}  The set $1,\theta',\theta'^2,\ldots, \theta'^{p-1}$ forms a basis of 
$S'/R'$. 
\end{lemma}

\begin{proof} It is enough to show that this set spans, 
and for that it is enough to show this set spans 
$S'$ modulo any maximal ideal, $M$, of $R'$. 
If $p \in M$ so $\eta \in M$ then 
$\sigma(\bar {\theta'}) = \bar {\theta'} + 1$ 
in $S'/MS'$. $S'/M'$ contains an element 
$\theta''$ such that $\sigma(\theta'') = 
\theta'' + 1$, we have $\bar {\theta'}= 
\theta'' + r$ 
for some $r \in R'$. Since the powers of $\theta''$ 
span it is clear that the powers of $\bar {\theta'}$ 
span also. 

If $p \notin M$ then $\eta$ is invertible in 
$R'/M$ and we set 
$\alpha = 1 + \eta\bar {\theta'}$. 
Then $\sigma(\alpha) = \rho\alpha$ 
and we compute that $\alpha = 
\sum_i \rho^{p-i+1}\sigma^i(z)$. In this case 
$S'/S'M$ is a Kummer extension of $R'/M$ 
and there is an $\alpha'$ such that 
$\sigma(\alpha') = \rho\alpha'$ and we can easily 
show that $\alpha' = r\alpha$ for some 
$r \in R'/M'$. Since the powers of $\alpha$ span 
it easily follows that the powers of $\alpha$ span. 
From this it easily follows that the powers of 
$\bar {\theta'} = (\alpha - 1)/\eta$ span.~\qed 

Returning to the proof of Theorem~\ref{thm:1.4}, 
let $a = (-1)^{p-1}\prod_i \sigma^i(\theta') \in R'$. 
It follows that $\theta'$ is a root 
of $Z^p + g(Z) - a$ and so $S' \cong 
R'[Z]/(Z^p + g(Z) - a)$. Since $S'/R'$ is 
separable, $Z^p + g(Z) - a$ has unit discriminant 
and $1 + a\eta^p$ is invertible. The theorem 
is clear.
\end{proof}

\begin{corollary}\label{cor:1.6}  Suppose $h: R' \to R''$ is a surjection of 
$Z[\rho]$ algebras such that $h^{-1}(1 + \eta^p{R''})$ are all units in 
$R'$. If $S''/R''$ is a $C_p$ Galois extension 
there is a $C_p$ Galois extension $S'/R'$ such that 
$S' \otimes_h R'' \cong S''$ as Galois extensions. 
\end{corollary}

\begin{proof} Since $R$ has the form $\Z[\rho][u](1/s)$ for 
$s \in 1 + \eta^p{R}$, any 
$\Z[\rho]$ algebra 
homomorphism $g'': R \to R''$ lifts to a $g': R \to R'$.
\end{proof}

\section{Doing without $\rho$}\label{sec:two} 

The question naturally arises whether we can perform a construction like that of the previous 
section without any assumption about $\rho$. 
The approach we use is much like the Kummer 
case, in that we adjoin a $\rho$, use the above 
section with the extra data, and thereby describe 
our $C_p$ Galois extensions. Let $\tau$ 
generate the Galois group of $\Q(\rho)$ 
thereby inducing an automorphism of $\Z[\rho]$. 
Let $s$ be such that $\tau(\rho) = \rho^s$. 
We can choose $s$ such that $s^{p-1} - 1 = 
pr$ and $r$ is prime to $p$. 

When $R'$ is an arbitrary commutative ring 
there are many ways to ``adjoin'' $\rho$ 
and it seems best to list the properties 
we need and then prove an existence theorem. 
We require that $R' \subset R$, that 
$R$ is a $\Z[\rho]$ algebra, that $\tau$ 
extends to $R$ of order $p-1$ or $1$, and 
that $R^{\tau} = R'$. Obviously if 
$pR' = 0$ we can set $R = R'$ and 
note that since $\Z[\rho]/(\eta) = 
\Z/p\Z$ we can view $R$ as an $\Z[\rho]$ 
algebra. However in any case we can 
define $R = \Z[\rho] \otimes_{\Z} R'$ 
and observe in the next result that $R$ 
has the properties we need. 

\begin{lemma}\label{lem:2.1}  The induced map of 
$\tau$ on the image of $\Z[\rho]$ in 
$R$ has order $p-1$. The fixed ring $R^{\tau} = R'$. 
\end{lemma}

\begin{proof} For the first statement, let 
$m\Z$ be the kernel of $\Z \to R'$. 
Then the image of $\Z[\rho]$ in $R$ is 
isomorphic to $\Z[\rho]/m\Z[\rho]$. 
As a group of prime order, the 
multiplicative group $<\rho>$ has image of 
order $p$ or 1. However $\eta = \rho - 1$ 
is never in $m\Z[\rho]$ for any $m$ so 
$<\rho>$ has image of order $p$. 
It follows that $\tau$ has order $p-1$ 
when acting on $\Z[\rho]/m\Z[\rho]$. 

As for the second statement, 
clearly $R = R'[Z]/(1 + Z + \cdots + Z^{p-1})$ 
and has basis $1,\rho,\ldots,\rho^{p-2}$ over $R'$, 
where $\rho$ is the image of $Z$. Also, 
$\tau(\rho) = \rho^s$. Let $\sum_0^{p-2} r_i\rho^i$ 
be an arbitrary element fixed by $\tau$ and hence any power 
of $\tau$. Assume $\tau^k(\rho) = \rho^t$. 
Then $\sum_i r_i\rho^i = \sum_i r_i\rho^{ti}$. 
We take the subscripts and exponents modulo $p$, 
and so it makes sense to write the right hand side 
as $\sum_i r_{i/t}\rho^i$ but note that $\rho^{-1}$ 
does not appear on the left side but $\rho^{-1} = 
-(\sum_i \rho^i)$. Thus the $\tau^k$ invariance can 
be written $\sum_0^{p-2} r_i\rho^i = 
\sum_0^{p-1} (r_{i/t} - r_{-1/t})\rho^i$. 
Looking at the $\rho^0$ terms we have $r_0 = 
r_0 - r_{-1/t}$ or $r_{-1/t} = 0$. Looking 
at all the possible values for $t$ we have 
$r_j = 0$ for $j > 0$ as needed.
\end{proof}

We begin with 
$S'/R'$ $C_p = <\sigma>$ Galois, and embed $R' \subset R$ a $\Z[\rho]$ algebra such that $\tau$ 
extends to $R$ with order $p-1$ or 1 as appropriate. If we set $S = S' \otimes_{R'} R$ 
then $S/R$ is also $<\sigma>$ Galois, 
$\tau$ extends (as $1 \otimes \tau$) to 
$S$ with the same order, and $\tau\sigma = 
\sigma\tau$. The question we next discuss 
is the reverse. If $R \supset R'$ is as above, 
and $S = R[Z]/(Z^p + g(Z) - u)$ is Galois 
with group $<\sigma>$, when can we extend 
$\tau$ to $S$ with the same order such that 
$\sigma\tau = \tau\sigma$? 

In the Kummer case this is known (\cite{S1}), but 
we recall it here for motivation. 
Given a Kummer element $\alpha$ with 
$\sigma(\alpha) = \rho\alpha$, 
then $\sigma\tau = \tau\sigma$ implies 
that $\sigma(\tau(\alpha)) = \tau(\sigma(\alpha)) = \tau(\rho\alpha) = \rho^s\tau(\alpha)$ and 
so $\tau(\alpha) = \alpha^s/z'$ for $z'$ $\sigma$ 
fixed. If we alter $\alpha$ and $\sigma$, 
we can assume $z' = z^r$. 
For any element $x$, define $M(x) = 
x^{s^{p-2}}\tau(x)^{s^{p-3}}\cdots\tau^{p-2}(x)$.   Using the fact 
that $\tau^{p-1} = 1$, we have $\alpha = 
\tau^{p-1}(\alpha) = \alpha^{s^{p-1}}/M(z)^r$ 
or $\alpha^{pr} = M(z)^r$. After a root of unity 
adjustment, $\alpha^p = M(z)$. Conversely, 
if $\alpha$ satisfies this equation we can 
define $\tau(\alpha) = \alpha^s/z^r$ and check 
that this defines an extension of $\tau$ such that $\tau\sigma = \sigma\tau$ and $\tau^{p-1} = 1$. 

We cannot use the above argument in our case because $\eta$ may be a zero divisor and 
so $\phi: R \to 1 + R\eta$ not injective. 
Instead assume $S = R[Z]/(Z^p + g(Z) - u)$ 
is Galois over $R$ 
with group generated by $\sigma$ 
with $\sigma(\theta) = \rho\theta + 1$. 
We wrote $1 + R\eta$ as $A$ above and set 
$A^*$ to be the units in $1 + R\eta$. 
in a similar way we define $B = 1 + S\eta$ 
and $B^*$, as well as $B_p = 1 + S\eta^p$ 
and $B_p^*$. 

\begin{lemma}\label{lem:2.2}  $B^{\sigma} = A$. 
\end{lemma}

\begin{proof} This is obvious if $\eta$ is a 
non--zero 
divisor so let $J \subset R$ to be the 
annihilator of $\eta$. If $J'$ is the 
annihilator of $\eta$ in $S$ then $J'$ 
is $\sigma$ invariant and since $S/R$ 
is Galois, $J' = (J' \cap R)S = JS$ and 
$S/J'$ is Galois over $R/J$ with group 
$<\sigma>$. In particular, $(S/J')^{\sigma} = 
R/J$. If $1 + x\eta$ is $\sigma$ fixed 
then $\sigma(x)\eta - x\eta = 0$ or 
$\sigma(x) = x + j$ where $j \in J'$. 
Thus $x$ maps to $(S/J')^{\sigma} = R/J$ 
and $x = r + j'$ where $j' \in J'$ and 
$r \in R$. But $1 + x\eta = 1 + r\eta$.
\end{proof}

Note that since $S/R$ is Galois we can conclude that $1 + u\eta^p$ 
and hence $1 + \theta\eta$ are units. 
If $\tau$ 
has order 1, this is just Artin-Schreier 
theory and so we will assume in this discussion 
that $\tau$ has order $p-1$. 
Since the image of $\rho$ has order $p$, 
we will drop the notation $\rho'$ 
and use $\rho$ to mean the element of 
$\Z[\rho]$ AND the image of $\rho$ in a 
$\Z[\rho]$ algebra $R$. 

Let $\delta = 
\delta_s = (\rho^s - 1)/(\rho - 1)$. 
Note that $\delta$ is a unit in $\Z[\rho]$ 
and hence in any $\Z[\rho]$ algebra. 
Also, $\tau(\eta) = \delta\eta$. 
If $\phi: S \to B$ is $\phi(z) = 1 + z\eta$, 
it follows that $\tau(\phi(z)) = 
\phi(\delta\tau(z))$. If $\phi_p: S \to 
B_p$ is $\phi_p(z) = 1 + z\eta^p$, then 
$\tau(\phi_p(z)) = \phi_p(\delta^p\tau(z))$. 
For the above reasons we sometimes will shift to studying the morphisms $\delta\tau$ 
and $\delta^p\tau$. Note that the 
$\tau$ norm $N_{\tau}(\delta) = 
(\rho^s - 1)/(\rho - 1)(\rho^{s^2} - 1)/(\rho^s - 1)\cdots(\rho^{s^{p-1}} - 1)/(\rho^{s^{p-2}} - 1) = 
1$ so $(\delta\tau)^{p-1} = 1$ and similarly 
$(\delta^p\tau)^{p-1} = 1$. 
Since $\delta \in \Z[\rho]$ if $\sigma$ and 
$\tau$ commute then $\sigma$ also commutes 
with $\delta\tau$ and $\delta^p\tau$.

We will also need to detail how $\tau$ interacts 
with the $\#_p$ operation from the last 
section. We collected the easy relations in the Appendix. 

Now assume $\tau$ and $\sigma$ act on $S$ as above. That is, $\tau\sigma = \sigma\tau$, 
$\tau^{p-1} = 1$, and $\tau$ extends the usual $\tau$ on $\Z[\rho]$. $S$ is generated by a $\theta$ such that 
$\sigma(\theta) = \rho\theta + 1 = 
\theta \oplus 1$. Then 
$\sigma(\delta(\tau(\theta))) = 
\delta\tau\sigma(\theta) = \delta\tau(\theta \oplus 1) = \delta(\tau(\theta)) \oplus 
\delta\tau(1) = \delta(\tau(\theta)) \oplus \delta$. However, we saw above that $\sigma(s*\theta) = 
s*\theta \oplus \delta$. Moreover, $\delta\tau(\theta) \in B^*$ so we can set 
$z = (s*\theta) \ominus \delta\tau(\theta)$ 
and compute that $\sigma(z) = z$. 
That is, there is a $z \in A^*$ 
such that $\delta\tau(\theta) = (s * \theta) \ominus z$. 

We will now engage in a series of arguments 
that will show that we can successively 
``improve'' $z$. 
Applying $\delta\tau$ again we 
have $(\delta\tau)^2(\theta) = (s^2 * \theta) 
\ominus (s * z)$. By induction we have 
$$\theta = (\delta\tau)^{p-1}(\theta) = 
(s^{p-1} * \theta) \ominus N'(z)$$ 
where $N'(z) = (s^{p-2} * z) \oplus 
(s^{p-3} * \delta\tau(z)) \oplus \cdots \oplus 
(\delta\tau)^{p-2}(z)$. Rewriting the set off 
equation we have $N'(z) = 
(s^{p-1} - 1) * \theta = p * (r * \theta)$. 
We saw above that $p * x \in B^*_p$ where 
$B^*_p$ are the units in $1 + \eta^pS$. 

Now $p * B^* \subset B^*_p$ so $B^*/B^*_p$ 
is a module over 
$F_p[<\delta\tau>]$ 
and this group ring is semisimple. Moreover, 
the irreducible modules of 
$F_p[<\delta\tau>]$ are 
all of the form $V_a = F_pv_a$ where $a \in F_p^*$
and $\delta\tau(v_a) = av$. Of course, 
$B^*/B^*_p$ is a direct sum of irreducibles. 
We compute $N'(v_a) = s^{p-2}v_a + s^{p-3}av_a + 
\cdots + a^{p-2}v_a =  (s^{p-2} + s^{p-3}a + \cdots + a^{p-2})v_a = 
a^{p-2}((b^{p-1} - 1)/(b - 1))v_a$ where 
$b = s/a$. Thus $N'(v_a) = 0$ except when 
$a = s$ and in this case 
$N'(v_s) = (p-1)s^{p-2} = -(1/s)v_a \not= 0$. 
Thus the fact that $N'(z) \in B^*_p$ implies 
that $z$ only has nonzero components in $V_a$ 
for $a \not= s$. However, if $v \in V_a$ 
for $a \not= s$ then $v = (\delta\tau - s)v' 
= (a - s)v'$ for some $v'$. Thus 
$z = (\delta\tau(y) - (s * y)) \oplus z'$ 
for some $z' \in B^*_p$. To state the next result we define $r *_p z$ via $\phi_p(z)^r = \phi_p(r *_p z)$ or $r *_p z = 
z \oplus_p \cdots \oplus_p z$ ($r$ times). 

\begin{lemma}\label{lem:2.3} By altering $\theta$ 
we may assume $(s * \theta) \ominus 
\delta\tau(\theta))\in r *_p B^*_p$. 
\end{lemma}

\begin{proof} Using the $y$ appearing above we compute $\delta\tau(\theta \ominus y) = 
(s*\theta) \ominus z \ominus \delta\tau(y) = 
s*(\theta \ominus y) \oplus (s*y) \ominus z 
\ominus \delta\theta(y) = s*(\theta \ominus y) \ominus z'$. Thus 
changing  $\theta$ to $\theta \ominus y$  
we may assume $z \in B_p^*$. 
Write $rr' = kp + 1$. We can replace 
$\theta$ by $\theta' = (kp + 1) * \theta$ because 
$kp * \theta$ is $\sigma$ fixed. 
But $(s * \theta') \ominus (\delta\tau(\theta'))$ 
is $rr'$ times $(s * \theta) \ominus 
\delta\tau(\theta))$ and $B_p^*$ is a subgroup of $B^*$.
\end{proof}

To sum up, our equation has now become 
(after changing $z$ again), 

$$\delta\tau(\theta) = 
(s * \theta) \ominus (r *_p z)\eta^{p-1} = 
(s * \theta) \ominus (r * (z\eta^{p-1})).$$

We will repeatedly apply $\delta\tau$ to the 
above equation. Since $r*_p$ appears above, 
and since $\tau(\phi_p(x)) = 
\tau(\delta^p\tau(x))$ it makes sense to define a new operator 
$$N(z) = (s^{p-2} *_p z) \oplus_p 
(s^{p-3} *_p \delta^p\tau(z)) \oplus_p \cdots 
\oplus_p (s *_p (\delta^p\tau)^{p-3}(z) \oplus_p (\delta^p\tau)^{p-2}(z).$$ 
Note that $\delta^p\tau(N(z)) = 
s *_p (N(z) \ominus_p (s^{p-2} *_p z)) \oplus z = s *_p (N(z)) \ominus (pr *_p z)$.  

Applying, as promised, $\delta\tau$ to the above equation (and making use of the identities in the appendix) we have 
$$(\delta\tau)^2(\theta) = \delta\tau(s * \theta) \ominus \delta\tau((r * (z\eta^{p-1})) = $$
$$s * (\delta\tau(\theta)) \ominus 
(r * (\delta^p\tau(z)\eta^{p-1}) = 
s * ((s * \theta) \ominus r * (z\eta^{p-1}) 
\ominus (r * (\delta^p\tau(z)\eta^{p-1})) = $$
$$(s^2 * \theta) \ominus (r * (s * z\eta^{p-1}) 
\ominus (r * (\delta^p\tau(z)\eta^{p-1}) =  
(s^2 * \theta \ominus 
(r * ((s *_p z) \oplus_p (\delta^p\tau(z))\eta^{p-1}).$$ 

Repeating this argument we have: 

$$\theta = (\delta\tau)^{p-1}(\theta) = 
(s^{p-1} * \theta) \ominus (r * (N(z)\eta^{p-1})).$$  
Thus, 
$$0 = (pr * \theta) \ominus  (r *_p N(z))\eta^{p-1} = $$
$$ = (pr \#_p \theta)\eta^{p-1} \ominus (r *_p N(z))\eta^{p-1} = 
((r *_p (p \#_p \theta)) \ominus_p (r *_p N(z))\eta^{p-1}$$ or 
$$\eta^{p-1}(r *_p (p \#_p \theta \ominus_p 
N(z) = 0.$$ 
Our next goal will be to change $\theta$ 
and $z$ to remove the $\eta^{p-1}$ and 
$r *_p$ factors from this relationship. 
To this end, set $n = (r *_p (p \#_p \theta)) \ominus_p (r *_p N(z))$ so $n\eta^{p-1} = 0$ and therefore $pn = 0$.  Set $d = (p \#_p \theta) \ominus N(z)$ so $n = r *_p d$. We compute how 
$\delta^p\tau$ acts on $d$ and therefore $n$. 
First $\delta^p\tau((p \#_p \theta)) =  
(p \#_p \delta\tau(\theta)) = $
$$(p \#_p (s * \theta \ominus (r *_p z)\eta^{p-1})) = $$
$$ = ((p \#_p (s * \theta)) \ominus_p 
(p \#_p (r *_p z)\eta^{p-1})) = $$
$$ = (s *_p (p \#_p \theta)) \ominus_p p *_p (r *_p z) = $$
$$ = s *_p ((r *_p (p \#_p \theta)) \ominus_p 
(pr *_p z).$$ 
Also, 
$\delta^p\tau(N(z)) =  s *_p (N(z)) \ominus_p 
(pr *_p z)$. 
The $pr *_p z$ terms cancel so 
$\delta^p\tau(d) = s *_p d$ and  

$$\delta^p\tau(n) = s *_p n.$$

If $m$ is any element of a $\Z[\rho]$ algebra 
such that 
$\eta^{p-1}m = pm = 0$ and 
$\delta^p\tau(m) = s *_p m$, then 
$N(m) = (p-1) *_p s^{p-2} *_p m  = 
(p-1)s^{p-2} *_p m$. Note $m \oplus_p m = 
m + m  + m^s\eta^p = 2m$, so more generally 
$t *_p m = tm$ 
so $(p-1)s^{p-2} *_p m = 
(p-1)s^{p-2}m = (-1/s)m$ as $pm = 0$. 
Thus if $m = (-s/r)n$, $r *_p N(m) = n$. 
Now changing $z$ to $z \oplus_p m$, 
we note that $(r *_p (z \oplus_p m))\eta^{p-1} = 
(r *_p z)\eta^{p-1}$ so it is still true that 
$\delta\tau(\theta) = (s * \theta) \ominus z\eta^{p-1}$ and we have $r *_p (p \#_p \theta) = r *_p N(z)$ 
or 

$$r *_p d = 0$$

where $d = (p \#_p \theta) \ominus N(z)$ is as above and $\phi_p(d)$ is a unit. Just as above $\delta^p\tau(d) = 
s *_p d$. If $rr' = 1 + pm$ then $d = 
pm *_p (-d) = p *_p e$ where $e = m *_p -d$. 
Note that $\delta^p\tau(e) = s *_p e$, $e$ is $\sigma$ fixed and $\delta\tau(e\eta^{p-1}) = \delta^p\tau\eta^{p-1} = (s *_p e)\eta^{p-1} = s * (e\eta^{p-1})$.  
Also $p *_p e = p \#_p (e\eta^{p-1})$ and so we set $\theta' = \theta \ominus (e\eta^{p-1})$ 
and compute that $p \#_p \theta' = N(z)$. 
Moreover, $\delta\tau(\theta') \ominus (s * \theta') = 
\delta\tau(\theta) \ominus \delta\tau(e\eta^{p-1}) \ominus (s * \theta) \oplus 
(s * (e\eta^{p-1})) = \delta\tau(\theta) \ominus (s * \theta)$. 

Thus we have shown: 

\begin{theorem}\label{thm:2.4}  Suppose $S/R$ is a $C_p = <\sigma>$ 
Galois extension where $R$ is a semilocal $\Z[\rho]$ 
algebra. Assume $\tau$ has image of order $p-1$ 
on the image of $\Z[\rho]$ and extends to $S$ with the 
same order. Assume $\sigma\tau = \tau\sigma$ on $S$. 
Then $S = R[Z]/(Z^P + g(Z) - u)$ where 
$u = N(z)$ for some $z \in R$ with $1 + z\eta^p$ 
a unit. There is a root of $Z^P + g(Z) - u$, $\theta$, 
such that $\delta\tau(\theta) = (s * \theta) \ominus z\eta^{p-1}$. 
\end{theorem}

To construct our generic extension, let $R' = 
\Z[\rho][x_0,\ldots,x_{p-2}]$ and define an extension 
of $\tau$ by setting $\tau(x_i) = x_i$. Let $z = 
\sum_i x_i\rho^i$ which we can think of as a generic 
element. Set $C$ to be the $\tau$ norm $N_{\tau}(1 + z\eta^p)$ 
and set $R = R'(1/C)$. For this $R$ we set 
$A = 1 + \eta{R}$, $A_p = 1 + \eta^p{R}$, 
and $A^*$, $A_p^*$ the invertible elements in these 
monoids. Since all the $\tau^i(1 + z\eta^p)$ 
are invertible it is true that $u = N(z)$ is in 
$A^*_p$. Form $S = R[Z]/(Z^p + g(Z) + u)$ and let 
$\theta$ be the image of $Z$ in $S$. Then $S/R$ 
is Galois with group $<\sigma>$ where $\sigma(\theta) = 
\rho\theta + 1 = \theta \oplus 1$. 
We set $\tau(\theta) = (1/\delta)((s * \theta) \ominus r *_p (z\eta^{p-1}))$. 

\begin{theorem}\label{thm:2.5}  The above defines an extension 
of the action of $\tau$ on to $S$ where $\tau^{p-1} = 1$ 
and $\sigma\tau = \tau\sigma$. Moreover, $S^{\tau}/R^{\tau}$ 
is Galois with group $<\sigma>$. 
\end{theorem}

\begin{proof} In this case $\eta$ is a non zero divisor so 
we note that if $\alpha = 1 + \theta\eta$ then the above 
is equivalent to $\tau(\alpha) = (1 + \theta\eta)^s/(1 + z\eta^p)$. 
Since $S/R$ is etale, $S$ is integrally closed and 
so $S/S'$ is an integral extension where $S' = R(\alpha)$. 
It follows that $\tau(S) = S$. 
The previously known arguments in the Kummer case 
show that $\tau\sigma = \sigma\tau$ and $\tau^{p-1} = 1$ 
on $S'$ and therefore on $S$. 
\end{proof}

Consider the ring extension 
$S^{\tau}/R^{\tau}$. Clearly the group 
$C = <\sigma>$ acts on $S^{\tau}$ with invariant 
ring $R^{\tau}$. The theorem is proven once we observe: 

\begin{proposition}\label{prop:2.6}  $S^{\tau}/R^{\tau}$ 
is Galois with group $C$. 
\end{proposition}

\begin{proof}  This is not immediate because $R/R^{\tau}$ 
is ramified and hence not etale. Note that 
$R$ is finitely generated as a module over 
$R^{\tau}$. In suffices 
to prove this proposition after localizing at all maximal 
ideals, $M$, of $R^{\tau}$. 
If $p \notin M$ then $R_M/R^{\tau}_M$ is etale 
and this is clear. If $p \in M$ then $R_M$ 
is semilocal and $\eta$ is in all maximal ideals. 

Consider $\epsilon = \theta + \tau(\theta) + 
\cdots + \tau^{p-2}(\theta) \in S^{\tau}$. Then 
$\sigma^i(\tau^j(\theta)) = 
\tau^j(\sigma^i(\theta)) = 
\tau^j(\rho^i\theta + \rho^{i-1} + \cdots + 1) = 
\rho^{s^ji}\tau^j(\theta) + \rho^{s^j(i-1)} 
+ \cdots + 1$. Modulo $\eta$ this is just 
$\tau^j(\theta) + i$. Thus, modulo $\eta$, 
$\sigma^i(\epsilon) = \epsilon + (p-1)i = 
\epsilon - i$. It follows that modulo 
the Jacobson radical of $R_M$, $\epsilon$ 
is an Artin-Schreier element and so 
$1,\epsilon,\epsilon^2,\ldots,\epsilon^{p-1}$ 
is a basis of $S_M/R_M$. Thus 
$S^{\tau} \otimes_{R^{\tau}} R \cong S$ 
and $S^{\tau}_M$ is free as a module over 
$R^{\tau}_M$. The natural map 
$\Delta(S^{\tau}/R^{\tau}) \to 
\End_{\R^{\tau}}(S^{\tau})$ has a determinant 
$d \in R^{\tau}_M$. This map is an isomorphism 
over $R$ and so $d$ is invertible in $R$ and 
hence $d$ is invertible in $R^{\tau}$.
We are done by \cite[page~81]{DI}.
\end{proof}

It is pretty clear that $R^{\tau} = \Z[x_0,\ldots,x_{p-2}](1/C)$. 
We are ready for: 

\begin{theorem}\label{thm:2.7}  Suppose $R_1$ is a semilocal commutative 
ring and $S_1/R_1$ is $C_p = <\sigma>$ Galois. Then there 
is a $\psi: R^{\tau} \to R_1$ such that $S_1 \cong 
S^{\tau} \otimes_{\psi} R_1$ as $C_p$ Galois extensions. 
\end{theorem}

\begin{proof} Form $R_1' \cong \Z[\rho] \otimes_{\Z} R_1$ 
and 
$S_1' = S_1 \otimes_{R_1} R_1' = \Z[\rho] \otimes_{\Z} S_1$. 
Then $\tau$ acts on $S_1'$ as $\tau \otimes 1$ 
and commutes with $\sigma$. We have $(R_1')^{\tau} = 
R_1$ and $(S_1')^{\tau} = S_1$. 

Since $R_1'/R_1$ is is finite, $R_1'$ is semilocal. 
In the nontrivial case $\tau$ 
has order $p-1$. Since $S_1/R_1$ is $C_p = <\sigma>]$ Galois 
there is a $\theta' \in S_1$ and a $z' \in R_1'$ 
such that $1 + z'\eta^{p-1}$ is invertible, 
$\sigma\tau(\theta') = (s * \theta') \ominus z'\eta^{p-1}$, 
and $\theta'^p + g(\theta') = N(z')$. Write $z' 
= \sum x_i'\rho^i$. Define $\psi: S \to S_1'$ be setting 
$\psi(x_i) = x_i'$, and $\psi(\theta) = \theta'$. 
The morphism $\psi$ clearly preserves $\tau$ and $\sigma$ 
actions and therefore defines $S^{\tau} \to S^{\tau}$ 
which restricts to $R^{\tau} \to R_1$. Since $S^{\tau}/R^{\tau}$ 
is Galois and etale, the rest is clear.
\end{proof}

Given Theorem~\ref{thm:2.7} we have the corresponding lifting property. 

\begin{corollary}\label{cor:2.8}  Let $g: R' \to R''$ be a surjection 
of commutative rings with $R''$ semilocal. 
Assume $z \in R'$ is a unit if $g(z)$ 
is a unit. 
Let $S''/R''$ be $C_p$ Galois. Then there is a $C_p$ 
Galois $S'/R'$ such that $S' \otimes_g R'' \cong S''$ as $C_p$ 
Galois extensions.  
\end{corollary}

The above result is easily extended to Galois groups of the 
form $C_p \oplus \cdots \oplus C_p$. 

\section{A curious algebra}\label{sec:three} 

Of course, with a new description of cyclic 
Galois extensions we get a new description 
of cyclic algebras. In detail, let $K$ be a field and $\rho \in K$ a primitive $p$ root of one unless $K$ has characteristic $p$ in 
which case $\rho = 1$. 
$L/K, \sigma$ be a $C_p$ Galois extension generated 
by $\alpha$ with $\sigma(\theta) = 
\rho\theta + 1$. Then the cyclic algebra 
$A = \Delta(L/K,\sigma,b)$ is generated by 
$\alpha$ and $\beta$ such that $\alpha^p + 
g(\alpha) = u$ for some $u$, $\beta^p = b$, and $\beta\alpha\beta^{-1} = 
\rho\alpha + 1$. The point of the above 
description is that it is characteristic free, though it assumes $\rho \in K$. 

In characteristic $p$ there is a different 
description of ``cyclic'' algebras in terms 
of deriviations. Now assume $K$ has characteristic $p$, so the above description 
of $A$ reduces to $\alpha^p - \alpha = u$, 
$\beta^p = b$ and $\beta\alpha\beta^{-1} = 
\alpha + 1$. This last relation can be rewritten as $\beta\alpha - \alpha\beta = 
\beta$. If we set $\gamma = \alpha\beta^{-1}$, then $\beta\gamma - \gamma\beta = 1$, 
and $\gamma^p = u/b$ as well as $\beta^p = b$ 
still. If $c = u/b$ then we have written 
$A$ as the differential crossed 
product $[c,b]$. 

The nice thing about these differential crossed product algebras is that they are 
Azumaya in complete generality. To state 
the results, let $R$ be an arbitrary commutative ring with $pR = 0$. 
Define the algebra $[c,b]$ as generated 
over $R$ by $\gamma,\beta$ subject 
to the relations $\gamma^p = c$, $\beta^p = b$, and $\beta\gamma - \gamma\beta = 1$. 
We have (\cite{KOS}): 

\begin{theorem}\label{thm:3.1}  $[c,b]$ is an Azumaya 
algebra of rank $p^2$ for any choice of 
$c,b$. Furthermore, any order $p$ element 
of the Brauer group $\Br(R)$ is a product 
(in the Brauer group) of classes determined 
by algebras of this form. 
\end{theorem}

Inspired by these differential crossed products, go back to the general $A$ 
generated by $\alpha$, $\beta$ subject to 
$\alpha^p + g(\alpha) = u$, $\beta^p = b$, 
and $\beta\alpha\beta^{-1} = 
\rho\alpha + 1$. The last relation can be rewritten 
as $\beta\alpha - \rho\alpha\beta = \beta$. 
As before substitute $\gamma = \alpha\beta^{-1}$. Now the last relation reads as $\beta\gamma - \rho\gamma\beta = 1$. 

Now we pull back and be as general as possible. Let $B$ be the following ring. 
Form the free ring $Z[\rho]<x,y>$ and let 
$I$ be the two sided ideal generated by 
$xy - \rho{yx} - 1$. Set 
$B = \Z[\rho]<x,y>/I$. 

\begin{proposition}\label{prop:3.2}  $x^p$ and 
$y^p$ are central in $B$. The center of 
$B$ is the polynomial ring 
$C = \Z[\rho][x^p,y^p]$. As a module 
over $C$, $B$ is free with generators 
$\{\beta^i\gamma^j | 1 \leq i,j \leq p-1\}$. 
\end{proposition}

\begin{proof} We will also use $x,y$ to denote 
the respective images in $B$. With easy variable changes, there is 
symmetry between $x$ and $y$ and for the 
first statement it suffices to show 
$x^p$ is central. An easy induction 
shows that $x^iy = \rho^iyx^i + (1 + \rho + \cdots + \rho^{i-1})x^{i-1}$ and the 
first statement follows. 

It is clear that $B$ is spanned by all the 
monomials $x^iy^j$. Form the 
field $K = \Q(\rho)(u,v)$ and let $B'$ 
be the cyclic algebra $\Delta(L/K, \sigma,v)$ 
generated by $\alpha$, $\beta$ 
satisfying $\alpha^p + g(\alpha) = u$, 
$\beta^p = v$, and $\beta\alpha\beta^{-1} = 
\rho\alpha + 1$. If $\gamma = \alpha\beta^{-1}$ we say above that $\beta\gamma - \rho\gamma\beta = 1$. 
Thus there is a ring homomorphism 
$\phi: B \to B'$ with $\phi(x) = \beta$ and 
$\phi(y) = \gamma$. Note that $\phi(x^p) 
= v$ and $\phi(y^p) = u/v = \gamma^p$ 
since $u/v$ is the norm of $\gamma$. 
It follows that $C= \Z[\rho][x^p,y^p]$ is a polynomial ring and that the $x^iy^j$ are a 
free basis for $B$ over $C$. Finally 
$C$ is the center of $B$ because this center 
is finite over $C$, is a subring of 
$\Q(\rho)(u,v)$, and $C$ is integrally closed.
\end{proof} 

Write $C = \Z[\rho][s,t]$ where $s = x^p$ 
and $y^p = t$. We are interested in the locus 
where $B/C$ is Azumaya. There is 
a standard homomorphism 
$\psi: B \otimes_C B^{\circ} \to \End_C(B)$ 
where $\psi(a \otimes b)(z) = azb$. 
Of course $B/C$ is Azumaya exactly where 
$\psi$ is an isomorphism. 
Since the domain and range of $\psi$ are 
free $C$ modules, there is an 
$f \in C$ which is its determinant. 
Thus $B/C$ is Azumaya exactly on the open 
set $D(f)$. 

We can make some deductions about $f$. 
Modulo $p$, $B/C$ is Azumaya because 
it becomes a differential crossed product. 
Thus the ideal $(f,p) = C$. 
Furthermore, by translating back to 
cyclic algebras, $B/C$ is Azumaya when 
$s$ and $1 + st\eta^p$ are invertible. 
Thus $f$ is a divisor of $s(1 + st\eta^p)$. 
But by symmetry $f$ is also a divisor 
of $t(1 + st\eta^p)$. Finally, 
$B/C$ is NOT globally Azumaya. 
If $I = (1 + st\eta^p)$, then  
the image of $1 + \eta{xy}$ generates a
nilpotent two sided ideal in $B/IB$. 

We have shown, 

\begin{theorem}\label{thm:3.3}  $B/C$ is Azumaya 
precisely on the open set $D(1 + st\eta^p)$. 
\end{theorem}

We next observe: 

\begin{theorem}\label{thm:3.4} $B/C$ 
is a maximal order. 
\end{theorem}

\begin{proof} By \cite{AG} Theorem 1.5, 
it suffices to prove that for 
any height one prime $P \subset C$, 
$B_P/C_P$ is a maximal order, 
where $C_P$ is the localization of 
$C$ at the prime $P$. 

We claim $1 + st\eta^p$ is a prime 
element, meaning that $(1 + st\eta^p)$ is a height one prime. However 
$\bar C = C/(1 + st\eta^p) \cong 
\Z[\rho][t](1/t\eta^p)$ is a 
domain making this clear. 

If $P$ is a height one prime that 
is not $(1 + st\eta^p)$, then 
$B_P/C_P$ is Azumaya and therefore 
a maximal order. If $P = (1 + st\eta^p)$, let $\hat C_P$ be the 
completion of the dvr $C_P$. 
It suffices to show that 
$\hat B_P = B \otimes_C \hat C_P$ 
is a maximal order. However, 
$\hat B_P$ is the cyclic algebra 
generated by $y,z$ where $z = 
1 + yx\eta$. Note that $y^p = t$, 
$z^p = 1 + st\eta^p$, and 
$y^{-1}zy = 1 + xy\eta = 
1 + (1 + \rho{yx})\eta = 1 + \eta + 
\rho(\eta{yx}) = \rho(1 + yx\eta)$. 
Note that this $\hat B_P$ has 
ramification at $P$ defined by 
$t^{1/p}$ which is a nontrivial 
extension of the fraction field, 
$\hat K$,  of 
$\bar C$. This implies $D = B \otimes_C \hat K$ is a division algebra with $\hat B_P$ as an order. Since $1 + st\eta^p$ is the prime element of 
$\hat C_P$, it follows that $z = 
1 + xy\eta$ is the prime element 
of $D$ and $\hat B_P$ is the standard maximal order.
\end{proof} 

Tautologically, $B(1/(1 + st\eta^p))$ 
is some sort of generic object, since it is generic for all its specializations. 
Instead of specifying all of this, let us 
note the obvious associated surjectivity 
result. 

\begin{corollary}\label{cor:3.5}  Suppose $R$ is an 
arbitrary $\Z[\rho]$ algebra with ideal $I$ such that 
$p \in I$ and all elements in 
$1 + pR$ are invertible. Then the induced map of Brauer groups $\Br(R) \to 
\Br(R/I)$ is surjective when restricted 
to elements of order $p$. 
\end{corollary}

\begin{proof} First note that as $(1 + z\eta)^p 
\in 1 + pR$ for any $z \in R$, we have that 
all elements of $1 +\eta{R}$ are inverible. 
To finish, it suffices to show that any 
degree $p$ differential crossed product 
over $R/I$ lifts to a degree $p$ Azumaya 
algebra over $R$, and this is obvious using 
$B$.
\end{proof}

One would like a variant of the algebra $B$, 
or just the surjectivity result 3.5, 
which does use or need $\rho$. The fact that 
$\Z[\rho]/\Z$ is not etale makes the construction of an explicit algebra difficult, but we can use corestriction on the field level to construct a Brauer class. 

Let $\tau$ extend to an action on $C$ 
by setting $\tau(s) = s$ and $\tau(t) = t$. 
Clearly $C^{\tau} = \Z[s,t]$ which we call $C'$. If $\tau^i(\rho) = \rho^r \not= \rho$
then set = $\pi = 1 + st\eta^p$ and 
compute $\tau^i(\pi) = 1 + st\delta_r^p\eta^p$ where we recall 
that $\delta_r = 1 + \rho + \ldots + \rho^{r-1}$. It follows that 
$\tau^i(\pi)$ and $\pi$ 
generate distinct prime ideals. 

Consider $\pi' = N_{\tau}(\pi) \in C'$. 
It is clear that $\pi'$ has the form 
$1 + wp$ for some $w \in C'$ and $\pi'$ is a prime element of $C'$ which 
splits completely in $C$. 

Now let $K = q(C)$ be the fraction field 
of $C$ which is Galois with group ${<\tau>}$ 
over $K' = q(C')$. Let $\beta \in \Br(K)$ 
correspond to $B \otimes_C K$ and form 
the corestriction $\alpha = \Cor_{K/K'}(\beta) \in \Br(K')$. We saw above that 
$\beta$ only ramifies at one height one 

prime namely $(\pi)$. It follows that 
$\alpha$ only ramifies at 
$(\pi')$. Because of the recent important result (\cite {KC}), 
we know the "purity of branch loci 
for Brauer groups" in this situation.  
Thus we have  
that $\alpha$ is the image of a Brauer 
class in $\Br(C'(1/\pi'))$ we also call 
$\alpha$. Since $C'(1/\pi')$ is regular, 
$\Br(C'(1/\pi')) \to \Br(K')$ is injective 
and $\alpha$ is unique. In a similar 
way we write $\beta \in \Br(C(1/\pi)$ 
to be the unique preimage of $\beta$. 

Let $R$ be a commutative ring with 
$pR = 0$. View $R$ as a $\Z[\rho]$ algebra 
via the morphism $\Z[\rho] \to 
\Z[\rho]/(\eta) = \Z/p\Z \subset R$. 
Let $\gamma \in \Br(R)$ be a Brauer class 
generated by a degree $p$ differential crossed product. We saw above that there 
is a $\Z[\rho]$ algebra homomorphism $f: C(1/\pi) \to R$ 
such that the induced Brauer group 
map $f: \Br(C(1/\pi)) \to \Br(R)$ 
satisfies $f(\beta) = \gamma$. 
Note that $f(\pi) = 1 = f(\tau^i(\pi))$. 
Thus $f$ extends to $f: C(1/\pi') \to R$. 
Since $f(\Z[\rho]) = \Z/p\Z$ it follows that 
$f \circ \tau^i = f$ for all $i$. Finally, 
$N_{\tau}(\beta)$ is the image of 
$\alpha$ in $\Br(C(1/\pi'))$. 
If we let $f': C'(1/\pi') \to R$ be the 
restriction of $f$, we can conclude that 
$f'(\alpha) = \gamma^{p-1} = \gamma^{-1}$. 
We now have: 

\begin{theorem} \label{thm:3.6} Suppose $R$ 
is a commutative ring and $I \subset R$ 
an ideal such that $p \in I$. Further assume 
that all elements of $1 + pR$ are units in $R$. Then $\Br(R) \to \Br(R/I)$ 
is surjective on elements of order $p$. 
\end{theorem}

\clearpage

\appendix

\section*{Appendix}

We gather here the rules relating $\oplus$, 
$\oplus_p$, $*$, $*_p$, $\#_p$, and how they interact 
with $\tau$. Recall that $\delta = \delta_s = 
(\rho^s - 1)/(\rho - 1)$. All the proofs are easily made 
by proving the relationship in the generic case, where 
$\eta$ is a non zero divisor, and then concluding the 
relation in general. 

$$\phi_p(x) = \phi(x\eta^{p-1}) \leqno{(1)}$$.

$$s*(x \oplus y) = (s*x) \oplus s*y) \leqno{(2)}$$ 
and 
$$s *_p( x \oplus_p y) = (s *_p x) \oplus_p 
(s *_p y). \leqno{(3)}$$ 

$$(pr \#_p x)\eta^{p-1} = (pr * x)\leqno{(4)}$$ 

and 

$$(p *_p x)= p \#_p (x\eta^{p-1}).\leqno{(5)}$$ 

$$((pr \#_p x) = p \#_p (r * x) = r *_p (p \#_p x).\leqno{(6)}$$

Also 

$$ pr \#_p (a \oplus b) = (pr \#_p a) \oplus_p 
(pr \#_p b)\leqno{(7)}$$ 

$$ pr \#_p (a \ominus b) = (pr \#_p a) \ominus_p 
(pr \#_p b) \leqno{(8)}$$ and 

$$pr \#_p (t * x) = t *_p (pr \#_p x). \leqno{(9)}$$  

$$pr \#_p (x\eta^{p-1}) = pr *_p x. \leqno{(10)}$$

$$(a \oplus_p b)\eta^{p-1} = a\eta^{p-1} \oplus 
b\eta^{p-1} \leqno{(11)}$$  

$$ (t *_p x)\eta^{p-1} = t * (x\eta^{p-1}).\leqno{(12)}$$

$$ \delta\tau(x \oplus y) = (\delta\tau(x)) 
\oplus (\delta\tau(y)) \leqno{(13)}$$ 

$$ \delta^p\tau(x \oplus_p y) = \delta^p\tau(x) \oplus_p 
\delta^p\tau(y) \leqno{(14)}$$

$$ \delta\tau(t * x) = t * \delta\tau(x) \leqno{(15)}$$ 

$$ (\delta^p\tau)(t *_p x) = t *_p (\delta^p\tau)(x) \leqno{(16)}$$ 

$$ \delta^p\tau(pr \#_p x) = pr \#_p \delta\tau(x) \leqno{(17)}$$

$$ \delta\tau(x\eta^{p-1}) = (\delta^p\tau(x))\eta^{p-1}. \leqno{(18)}$$

\bigskip


\end{document}